\documentclass[11pt,a4paper]{article}
\usepackage{amsmath}
\usepackage{amsthm}
\usepackage{amsfonts}
\usepackage{amssymb}
\usepackage{stmaryrd}
\usepackage{latexsym}
\usepackage[usenames]{color}

\newtheorem{theorem}{Theorem}

\newtheorem{algorithm}{Algorithm}
\newtheorem{procedure}{Procedure}

\theoremstyle{definition}
\newtheorem{definition}{Definition}

\newtheorem{remark}{Remark}
\newtheorem{example}{Example} 

\begin{document}

\title{\Large{\textbf{On  some quasigroup cryptographical primitives}}}
\author{\normalsize {Piroska ~Cs\" org\H o$^1$ and  Victor Shcherbacov$^2$}}

 \maketitle

\begin{abstract}
\noindent Using Vojvoda approach \cite{VOIVODA_04} we demonstrate that cryptographical primitives proposed in \cite{Petrescu_10} are vulnerable  relative to chosen ciphertext attack  and chosen plaintext attack. We develop proposed in \cite{SCERB_09_CSJM} modifications and add some new modifications of known quasigroup based stream ciphers \cite{MARKOVSKI, Petrescu_10}.  Systems of orthogonal $n$-ary groupoids are used.

\medskip


\medskip

\end{abstract}


\bigskip

\section{Introduction}

\subsection{Basic definitions}

Information on quasigroups and $n$-ary quasigroups it is possible to find in \cite{VD, 2, 1a, HOP}, on ciphers in \cite{MENEZES}.
Stream-ciphers based on quasigroups and their parastrophes were  discovered in the end of the XX-th century \cite{MARKOVSKI}. See also \cite{OS, Grosek_Sys, SCERB_09_CSJM}.
We give some definitions.

A sequence $x_m, x_{m+1}, \dots, x_n$, where $ m, n$ are natural numbers and $m\leq n$,  will be denoted by
$x_m^n$. If $m > n$,  then $x_m^n$ will be considered empty.  The expression $\overline{1,
n}$ designates the set $\{1, 2, \dots, n\}$ of natural numbers \cite{2}.

A non-empty set $Q$  together with an $n$-ary operation $A : Q^n \rightarrow Q$, $n\geq 2$
is called $n$-groupoid and it is denoted by $(Q, A)$.

It is convenient to define $n$-ary quasigroup in the following manner.
\begin{definition} \label{def2} An $n$-ary groupoid $(Q, A)$ with $n$-ary operation $A$ such
that in the equality $A(x_1,$ $ x_2, \dots, x_n) = x_{n+1}$ the knowledge of any $n$  elements from the elements $x_1, x_2, \dots,$
$ x_n, x_{n+1}$ uniquely specifies the remaining one is called $n$-ary quasigroup  \cite{2}.
\end{definition}

We give a classical equational definition of binary quasigroup \cite{EVANS_49}.

\begin{definition}  \label{EQUT_QUAS_DEF}
 A binary groupoid $(Q, A)$ is called a binary quasigroup if on the set $Q$ there exist
operations ${}^{(13)}A$ and ${}^{(23)}A$ such that in the algebra $(Q, A, {}^{(13)}A, {}^{(23)}A)$ the following
identities are fulfilled:
\begin{equation}
A({}^{(13)}A(x, y), y) = x, \label{(2e)}
\end{equation}
\begin{equation}
{}^{(13)}A(A(x, y), y) = x, \label{(4e)}
\end{equation}
\begin{equation}
A(x, {}^{(23)}A (x, y)) = y, \label{(1e)}
\end{equation}
\begin{equation}
{}^{(23)}A(x, A (x, y)) = y. \label{(3e)}
\end{equation}
\end{definition}

By tradition the operation $A$ is denoted by $\cdot$, ${}^{(23)}A$ by $\backslash$ and ${}^{(13)}A$ by $\slash$.

It is well known \cite{VD, HOP, SCERB_03} that any quasigroup $(Q, A)$ defines else five quasigroups, so-called parastrophes of quasigroup $(Q,A)$, namely  $(Q,{}^{(13)}A)$, $(Q,{}^{(23)}A)$, $(Q,{}^{(12)}A)$, $(Q,{}^{(123)}A)$, $(Q,{}^{(132)}A)$.

It is possible to give equational definition of  $n$-ary quasigroup as a generalization of Definition  \ref{EQUT_QUAS_DEF}. We follow \cite{2, Petrescu_10}.

\begin{definition}
 An $n$-ary groupoid $(Q, A)$ is called an $n$-ary quasigroup if on the set $Q$ there exist
operations ${}^{(1,\, n+1)}A$,  ${}^{(2,\, n+1)}A$, $\dots$, ${}^{(n,\, n+1)}A$ such that in the algebra $(Q, A, {}^{(1,\,n+1)}A, \dots, {}^{(n,\,n+1)}A)$  the following identities are fulfilled for all $i \in \overline{1,n}$:
\begin{equation}
A(x_1^{i-1}, {}^{(i,\, n+1)}A(x_1^n), x_{i+1}^n) = x_i, \label{(1ne)}
\end{equation}
\begin{equation}
{}^{(i,\, n+1)}A(x_1^{i-1}, A(x_1^n), x_{i+1}^n) = x_i. \label{(2ne)}
\end{equation}
\end{definition}

\subsection{Quasigroup based cryptosystem} \label{Quasigroup based cryptosystem}

We give based on binary quasigroup  algorithm  for secure  encoding. We use \cite{SCERB_03}.

A quasigroup $(Q,\cdot)$ and its $(23)$-parastrophe $(Q,\backslash)$  satisfy the following identities   $ x\cdot (x\backslash y) = y$,
$ x\backslash (x\cdot y) = y$. This is identities (\ref{(1e)}) and (\ref{(3e)}), respectively.

The authors \cite{MARKOVSKI}
propose to use this quasigroup  property  to construct a stream cipher.

\begin{algorithm} \label{ALG1} Let $Q$ be a non-empty alphabet,  $k$ be a natural number, $u_i,
v_i \in Q$, $i\in \{1,..., k\}$.  Define a quasigroup $(Q, A)$.
It is clear that the quasigroup $(Q, {}^{(23)}A)$ is defined in a unique way.

Take a fixed element $l$ ($l\in Q$), which  is called a leader.

Let $u_1 u_2... u_k$ be a $k$-tuple  of letters from $Q$.

It is proposed  the following ciphering procedure   $v_1 = A(l, u_1) , v_{i}= A(v_{i-1}, u_{i})$, $i= 2,..., k$.

  Therefore we obtain the following cipher-text $v_1v_2 \dots v_k$.

  The enciphering algorithm is constructed in the following way: $u_1= {}^{(23)}A(l, v_1)$,
$u_{i}= {}^{(23)}A(v_{i-1}, v_{i}),$ $i = 2,..., k.$

Indeed ${}^{(23)}A(v_{i-1}, v_{i}) = {}^{(23)}A(v_{i-1}, A(v_{i-1}, u_{i})) \overset{(\ref{(3e)})}{=} u_i.$
\end{algorithm}

\begin{remark}
The equality $A = {}^{(23)}A$ is fulfilled if and only if  $A(x, A(x,y)) =y$ for all $x, y \in Q$.
\end{remark}

\begin{example}

 Let alphabet $A$ consists from the letters $a, b, c$. Take the quasigroup $(A, \cdot)$:
\[
{\begin{array}{c|ccc}
\cdot & a & b & c  \\
\hline
a & b & c & a   \\
b & c & a & b   \\
c & a & b & c
\end{array}}
\]
Then  $(A, \backslash)$ has the following Cayley table
\[
{\begin{array}{c|ccc}
\backslash & a & b & c  \\
\hline
a & c & a & b   \\
b & b & c & a   \\
c & a & b & c
\end{array}}
\]

Let $l=a$ and open text is $u = b\,b\,c\,a\,a\,c\,b\,a$. Then the cipher text  is $v = c\,b\,b\,c\,a\,a\,c\,a$.
Applying of decoding function on $v$ we get $b\,b\,c\,a\,a\,c\,b\,a=u$.
\end{example}

In \cite{OS} the authors claimed  that this cipher is resistant  to the brute force attack (exhaustive search) and to  the
statistical attack (in many languages some letters meet more frequently, than other letters).

In dissertation of Milan Vojvoda \cite{VOIVODA_04} it is proved that this cipher is not resistant to chosen ciphertext attack  and chosen plaintext attack. It is claimed that this cipher is not resistant relatively  statistical attack (Slovak language).

There exist few  ways to generalize Algorithm \ref{ALG1}.
The most obvious way is to increase arity of a quasigroup, i.e. instead of binary to apply $n$-ary ($n\geq 3$) quasigroups. This way was proposed in \cite{SCERB_03, SCERB_03_1} and  was realized in \cite{Petrescu_07, Petrescu_10}.
Notice Prof. A.~Petrescu writes that he found this $n$-ary generalization independently.

The second way was proposed in fact in \cite{SCERB_09_CSJM}. Namely instead of pair of binary quasigroups it  was proposed to use a system of $n$ $n$-ary orthogonal operations (groupoids).

\section{Cryptanalysis of  $n$-ary quasigroup cipher}

\begin{algorithm} \label{ALGn_ar} Let $Q$ be a non-empty alphabet,  $k$ be a natural number, $u_i,
v_i \in Q$, $i\in \{1,..., k\}$.  Define an $n$-ary quasigroup $(Q, A)$.
It is clear that any quasigroup $(Q, {}^{(i,\, n+1)}A)$ for any fixed value $i$ is defined in a unique way. Below for simplicity we put $i=n$.

Take fixed elements $l_1^{(n-1)(n-1)}$ ($l_i\in Q$), which  are called  leaders.

Let $u_1 u_2... u_k$ be a $k$-tuple  of letters from $Q$.

It is proposed  the following ciphering (encryption) procedure
\begin{equation} \label{seven}
\begin{split}
& v_1 = A(l_1^{n-1}, u_1),  \\
& v_2 = A(l_n^{2n-2}, u_2), \\
& \dots , \\
& v_{n-1} = A(l_{n^2 - 3n +3}^{(n-1)(n-1)}, u_{n-1}), \\
& v_{n}= A(v_{1}^{n-1}, u_{n}), \\
& v_{n+1}= A(v_{2}^{n}, u_{n+1}), \\
& v_{n+2}= A(v_{3}^{n+1}, u_{n+2}), \\
& \dots
\end{split}
\end{equation}
Therefore we obtain the following cipher-text $v_1v_2 \dots, v_{n-1}, v_{n}, v_{n+1}, \dots$.

  The deciphering algorithm also  is constructed similarly with binary case:
  \begin{equation} \label{eight}
\begin{split}
& u_1= {}^{(n,\, n+1)}A(l_1^{n-1}, v_1),  \\
& u_{2}= {}^{(n,\, n+1)}A(l_n^{2n-2}, v_{2}), \\
& \dots , \\
& u_{n-1} = {}^{(n,\, n+1)}A(l_{n^2 - 3n +3}^{(n-1)(n-1)}, v_{n-1})\\
& u_{n}= {}^{(n,\, n+1)}A(v_{1}^{n-1}, v_{n}),\\
& u_{n+1}= {}^{(n,\, n+1)}A(v_{2}^{n}, v_{n+1}), \\
& u_{n+2}= {}^{(n,\, n+1)}A(v_{3}^{n+1}, v_{n+2}), \\
& \dots
\end{split}
\end{equation}
Indeed, for example,  ${}^{(n,\, n+1)}A(v_{1}^{n-1}, v_{n}) = {}^{(n,\, n+1)}A(v_{1}^{n-1}, A(v_{1}^{n-1}, u_{n})) \overset{(\ref{(2ne)})}{=} u_n.$
\end{algorithm}

Probably there exists a sense to use in Algorithm  \ref{ALGn_ar}  irreducible 3-ary or 4-ary finite quasigroup \cite{2, BORISENKO_V_V, AKIV_GOLD_00, AKIV_GOLD_01}.

\begin{remark}
In equation (\ref{seven}) (encryption procedure) and, therefore, in decryption procedure (equation (\ref{eight})) it is possible to use more than one  $n$-quasigroup operation.
\end{remark}

\subsection{Chosen ciphertext attack}

We describe chosen ciphertext attack on cipher defined in Algorithm \ref{ALGn_ar}. Binary analog of this attack  is described in \cite{VOIVODA_04}.
Let $Q = \{q_1, q_2, \dots , q_k\}$, $|Q|=q$ and assume the cryptanalyst has access
to the decryption device loaded with an unknown key. Then he can construct the following ciphertext:

$v_1 = q_1, v_2=q_1, \dots, v_{n-1} = q_1, v_{n} = q_1, v_{n+1} = q_2$. Then $u_n ={}^{(n,\, n+1)}A(q_1, \dots, q_1)$. If $v_{n+1} = q_2$,  then
$u_{n+1} ={}^{(n,\, n+1)}A(q_1, \dots, q_1, q_2)$, if $v_{n+2} = q_3$, then $u_{n+2} ={}^{(n,\, n+1)}A(q_1, \dots, q_1, q_2, q_3)$ and so on. Continuing in such manner we can find multiplication table of quasigroup $(Q, {}^{(n,\, n+1)}A)$, therefore multiplication table of quasigroup $(Q, A)$ too. Notice, $|(Q, A)| = q^n$.

Having multiplication table of quasigroup $(Q, {}^{(n,\, n+1)}A)$  we can easy encipher any ciphertext  starting from the symbol $v_n$.

\subsection{Chosen ciphertext attack on the  leader elements}

In order to decrypt the elements $v_1, \dots, v_{n-1}$ we should know action of $(n-1)$-tuples of  leader elements  on any element of the set $Q$.
In other words we should know the action of translation $T_1(l_1, l_2, l_{n-1}, -)$, $T_2(l_n, l_{n+1}, l_{2n-2}, -)$, $\dots$, $T_{n-1}$ on the set $Q$.

It is not difficult to find element-leader using quasigroup $(Q, {}^{(n,\, n+1)}A)$ in binary case. It is sufficient to solve  equation ${}^{(23)}A(l, a) = b$ for fixed elements $a, b \in Q$.

Notice for cryptographical purposes it is not necessary to find elements-leaders $l_1, l_2$. It is sufficient to find  pair of elements $c, d$ such that ${}^{(3\, 4)}A(l_1, l_2, x) = {}^{(3\, 4)}A(c, d, x)$ for all $x\in Q$. For this aims there exists a possibility to decrypt  $q$ one letter  cipher-texts $q_1, q_2, \dots, q_n$ in any order.

In order to establish the action of elements-leaders $l_3, l_4$ on the set $Q$ (action of translation $T(l_3, l_4, -)$) it is possible to decrypt $q$ pair of elements of the form  $a, q_1$, $a, q_2$, $\dots$, $a, q_q$, where $a$ is a fixed element of the set $Q$.

 It is possible to unite calculation of action of translations $T(l_1, l_2, -)$ and $T(l_3, l_4, -)$ in one procedure using by decryption $q$ pairs of elements $(q_i, q_j)$, where $q_i\neq q_j$, $\cup_{i=1}^{q}q_i = Q$, $\cup_{j=1}^{q}q_j = Q$.

In the similar way it is possible to operate in  $n$-ary case ($n\geq 4$).

\begin{example}\label{TERN_NO_REDUC_QUS}
We give an example of ternary quasigroup $(Q, A)$ of order 4 \cite[p.115]{2}. In some sense this quasigroup is non-trivial since it is not an isotope of $3$-ary group $(Q, f)$ with the form $f(x_1^3) = x_1+x_2+x_3$ where $(Q, +)$ is a binary group of order 4. Recall there exist two groups of order 4, namely cyclic group $Z_4$  and Klein group $Z_2\times Z_2$. Any binary quasigroup of order 4 is a group isotope \cite{A1,A2}.

\[
\begin{array}{cccc}
\begin{array}{c|cccc}
A_0 & 0 & 1 & 2 & 3 \\
\hline
0  & 0 & 1 & 2 & 3 \\
1  & 1 & 2 & 3 & 0 \\
2  & 2 & 3 & 0 & 1 \\
3  & 3 & 0 & 1 & 2 \\
\end{array}
&
\begin{array}{c|cccc}
A_1 & 0 & 1 & 2 & 3 \\
\hline
0  & 1 & 0 & 3 & 2 \\
1  & 0 & 1 & 2 & 3 \\
2  & 3 & 2 & 1 & 0 \\
3  & 2 & 3 & 0 & 1 \\
\end{array}
&
\begin{array}{c|cccc}
A_2 & 0 & 1 & 2 & 3 \\
\hline
0  & 2 & 3 & 0 & 1 \\
1  & 3 & 0 & 1 & 2 \\
2  & 0 & 1 & 2 & 3 \\
3  & 1 & 2 & 3 & 0 \\
\end{array}
&
\begin{array}{c|cccc}
A_3 & 0 & 1 & 2 & 3 \\
\hline
0  & 3 & 2 & 1 & 0 \\
1  & 2 & 3 & 0 & 1 \\
2  & 1 & 0 & 3 & 2 \\
3  & 0 & 1 & 2 & 3 \\
\end{array}
\end{array}
\]
Notice $A(0, 1, 2) = A_0(1, 2) = 3,$ $A(2, 3, 2) = A_2(3, 2) = 3.$ Moreover $A(0, 1, x) = A(2, 3, x)$ for any  $x\in Q$. Then translations $T(0, 1, -)$ and $T(2, 3, -)$ are equal, pairs of leaders $(0, 1)$ and $(2, 3)$ are equal from cryptographical point of view.
\end{example}

\subsection{Chosen plaintext attack}

Chosen plaintext attack is similar with  chosen ciphertext attack.

Let assume the cryptanalyst has access
to the encryption device loaded with an unknown key. Then he can construct the following plaintexts:

$u_1 = q_1, u_2=q_1, \dots, u_{n-1} = q_1, u_{n} = q_1, u_{n+1} = q_2$. Then $v_n =A(q_1, \dots, q_1)$. If $u_{n+1} = q_2$,  then
$v_{n+1} =A(q_1, \dots, q_1, q_2)$. If $u_{n+2} = q_3$, then $v_{n+2} =A(q_1, \dots, q_1, q_2, q_3)$ and so on. Continuing it in such
manner  we can find multiplication table of quasigroup $(Q, A)$. Notice $|(Q, A)| = q^n$.

\subsection{Chosen plaintext attack on the  leader elements}

Chosen plain text attack on   leader elements is similar with chosen ciphertext attack on   leader elements and we omit it.

\section{Ciphers based on orthogonal $n$-ary groupoids}

\subsection{Some definitions}

We give classical definition of orthogonality of  $n$-ary operations \cite{YAC, BEL_YAK}.

\begin{definition} \label{n_ORTHOGON}
 $n$-ary groupoids  $(Q, f_1)$, $(Q, f_2)$, $\dots$, $(Q, f_n)$ are called orthogonal, if for any fixed
$n$-tuple $a_1, a_2, \dots, a_n$ the following system of equations

$$\left\{ \begin{array}{l}
f_1(x_1, x_2, \dots , x_n) = a_1\\
f_2(x_1, x_2, \dots , x_n) = a_2 \\
\dots \\
f_n(x_1, x_2, \dots , x_n) = a_n
 \end{array} \right.
 $$
 has a unique solution.
\end{definition}

There exist various generalizations of definition of orthogonality of $n$-ary operations.

\begin{definition} \label{K_ORTHOG}
 $n$-ary groupoids  $(Q, f_1)$, $(Q, f_2)$, $\dots$, $(Q, f_k)$ ($2 \leq k\leq n$)  given on a set $Q$ of order $m$ are called orthogonal if the system of equations \[\left\{
\begin{array} {l}
f_1(x_1, x_2, x_n)= a_1 \\
f_2(x_1, x_2, x_n)= a_2 \\
\dots \\
f_k(x_1, x_2, x_n)= a_k \\
\end{array}
\right. \]
 has exactly $m^{n-k}$ solutions for any k-tuple  $a_1, a_2, \dots, a_k$, where $a_1, a_2, \dots, a_k \in Q$ (see \cite{Bel_05}).
\end{definition}

If $k=n$, then from Definition \ref{K_ORTHOG} we obtain standard Definition \ref{n_ORTHOGON}.
Definition of orthogonality of binary systems has rich and  long history \cite{DK1}. About $n$-ary case see,
for example, \cite{EVANS_1}.

\begin{example} \label{SYST_EQ_3}
Operations $A_1(x_1, x_2, x_3) = 1\cdot x_1 + 0\cdot x_2 + 0\cdot x_3$, $A_2(x_1, x_2, $ $x_3) = 0\cdot x_1 +
1\cdot x_2 + 0\cdot x_3$, $A_3(x_1, x_2, x_3) = 0\cdot x_1 + 0\cdot x_2 + 1\cdot x_3$ defined over the  field $R$ of real
numbers (or over a finite field) are orthogonal, since the system
$$\left\{ \begin{array}{l}
1\cdot x_1 + 0\cdot  x_2 + 0\cdot x_3 = a_1\\
0 \cdot x_1 + 1\cdot x_2 + 0\cdot x_3 = a_2 \\
0\cdot x_1 + 0\cdot x_2 + 1\cdot x_3 = a_3
 \end{array} \right.
 $$
has a unique solution for any fixed $3$-tuple $(a_1, a_2, a_3)\in R^3$.
\end{example}

Notice any pair of ternary operations from Example \ref{SYST_EQ_3} is orthogonal in sense of Definition \ref{K_ORTHOG}.

We follow ideas of V.D. Belousov \cite{SYSTEMS_ORTHOGON}. See also \cite{BEL_YAK, YAC}.
It is easy to see that any system of $n$ orthogonal $n$-ary groupoids $(Q, f_i)$ $i\in \overline{1, n}$, defines a permutation of the set $Q^n$ and vice versa.  Thus there exist $(q^n)!$  $n$-ary orthogonal systems on a set   of order $q$.

\subsection{Construction of orthogonal $n$-ary groupoids}

In the following example will be given a sufficiently convenient and  general  way for the construction of systems of orthogonal $n$-ary groupoids.

\begin{example}
  Define operations $A_1(x_1, x_2, x_3)$, $A_2(x_1, x_2, x_3)$,  $A_3(x_1, x_2, x_3)$  over the set $M = \{ 0, \,  1, \, 2 \, \}$  in the following way. Take all $27$ triplets $ K = \{ (R_i,\, S_i, \, T_i) \, \mid \, R_i, S_i, T_i \in M,  i \in \overline{1, 27}\}$ in any fixed order and put
\begin{equation*}
\begin{split}
& A_1(0, 0, 0) = R_1, A_1(0, 0, 1) = R_2, A_1(0, 0, 2) = R_3, \dots,  A_1(2, 2, 2) = R_{27}, \\
& A_2(0, 0, 0) = S_1, A_2(0, 0, 1) = S_2, A_2(0, 0, 2) = S_3, \dots,  A_2(2, 2, 2) = S_{27}, \\
& A_3(0, 0, 0) = T_1, A_3(0, 0, 1) = T_2, A_3(0, 0, 2) = T_3, \dots,  A_3(2, 2, 2) = T_{27}.
\end{split}
\end{equation*}
The operations $A_1$, $A_2$ and $A_3$ form a   system of orthogonal operations. If we take this $27$ triplets in other order, then we obtain other system of orthogonal $3$-ary groupoids.
\end{example}

This way gives a possibility to construct easy inverse system $B$ of orthogonal $n$-ary operations to a fixed system $A$ of orthogonal $n$-ary operations. Recall inverse system means that $B(A(x_1^n)) = x_1^n$, $x_i \in Q$.

\begin{example}
We give example of three orthogonal ternary groupoids that are defined on four-element set $\{0, 1, 2, 3 \}$. Multiplication table of the first groupoid (in fact, of a quasigroup) is given in Example \ref{TERN_NO_REDUC_QUS}. Below we give multiplication tables of other two $3$-ary groupoids.

\[
\begin{array}{cccc}
\begin{array}{c|cccc}
A_0 & 0 & 1 & 2 & 3 \\
\hline
0  & 3 & 0 & 1 & 3 \\
1  & 0 & 2 & 3 & 0 \\
2  & 1 & 2 & 1 & 3 \\
3  & 1 & 1 & 2 & 2 \\
\end{array}
&
\begin{array}{c|cccc}
A_1 & 0 & 1 & 2 & 3 \\
\hline
0  & 2 & 1 & 1 & 0 \\
1  & 2 & 3 & 3 & 0 \\
2  & 0 & 2 & 1 & 3 \\
3  & 0 & 0 & 3 & 1 \\
\end{array}
&
\begin{array}{c|cccc}
A_2 & 0 & 1 & 2 & 3 \\
\hline
0  & 1 & 2 & 0 & 0 \\
1  & 2 & 0 & 3 & 1 \\
2  & 0 & 2 & 3 & 2 \\
3  & 3 & 2 & 1 & 1 \\
\end{array}
&
\begin{array}{c|cccc}
A_3 & 0 & 1 & 2 & 3 \\
\hline
0  & 3 & 3 & 2 & 2 \\
1  & 0 & 1 & 2 & 1 \\
2  & 0 & 2 & 0 & 3 \\
3  & 3 & 1 & 0 & 3 \\
\end{array}
\end{array}
\]
\[
\begin{array}{cccc}
\begin{array}{c|cccc}
A_0 & 0 & 1 & 2 & 3 \\
\hline
0  & 3 & 1 & 2 & 0 \\
1  & 2 & 1 & 1 & 2 \\
2  & 0 & 1 & 0 & 1 \\
3  & 3 & 1 & 2 & 3 \\
\end{array}
&
\begin{array}{c|cccc}
A_1 & 0 & 1 & 2 & 3 \\
\hline
0  & 1 & 2 & 1 & 3 \\
1  & 1 & 2 & 3 & 1 \\
2  & 0 & 2 & 2 & 0 \\
3  & 1 & 3 & 1 & 1 \\
\end{array}
&
\begin{array}{c|cccc}
A_2 & 0 & 1 & 2 & 3 \\
\hline
0  & 3 & 3 & 0 & 0 \\
1  & 2 & 1 & 0 & 1 \\
2  & 3 & 3 & 2 & 0 \\
3  & 3 & 0 & 2 & 3 \\
\end{array}
&
\begin{array}{c|cccc}
A_3 & 0 & 1 & 2 & 3 \\
\hline
0  & 2 & 1 & 0 & 0 \\
1  & 2 & 0 & 2 & 3 \\
2  & 3 & 3 & 2 & 0 \\
3  & 2 & 0 & 0 & 3 \\
\end{array}
\end{array}
\]
\end{example}

From formula $(q^n)!$ it follows that there exist $(4^3)! = 64!$   orthogonal systems of $3$-ary groupoids over a set of order 4.

\subsection{Ciphers on base of orthogonal systems of $n$-ary operation}

Here we propose to use a system  of orthogonal $n$-ary  groupoids as additional procedure in order to construct almost-stream cipher \cite{SCERB_09_CSJM}.

 Orthogonal systems of $n$-ary quasigroups were studied in \cite{Stojakovic_86, Syrbu_90, DUD_SYRB}.
Such  systems  have  more uniform distribution of elements of base set and therefore such systems may be more preferable in protection against statistical cryptanalytic  attacks.

\begin{procedure} \label{ALG3} Let $A$ be a non-empty alphabet,  $k$ be a natural number, $u_i,
v_i \in A$, $i\in \{1,..., k\}$.

\begin{enumerate}
  \item Define a system  of $n$ n-ary orthogonal operations $(A,f_i)$, $i = 1, 2, \dots, n$.
  \item We propose  the following enciphering procedure   $v_i = f_i(u_1, u_2, \dots, u_n)$, $i = 1, 2,..., n$. If $k<n$, then we can repeat plaintext or a part of plaintext necessary number of times.
\item  It is possible to apply the ciphering procedure more than one time. Number of applications of Step 2 can be non-fixed.
    \item Therefore we obtain the following ciphertext $v_1v_2 ...v_n$.
\end{enumerate}

The deciphering algorithm is based on the fact that  orthogonal system of n n-ary operations
\[ \left\{
  \begin{array} {l}
    f_1(x_1, x_2, \dots , x_n) = a_1 \\
  f_2(x_1, x_2, \dots , x_n) = a_2 \\
  \dots \\
    f_n(x_1, x_2, \dots , x_n) = a_n
    \end{array}
\right.
\]
  has a unique solution for any tuple of elements $a_1, \dots , a_n$.
 \end{procedure}

\subsection{Modifications of Procedure \ref{ALG3}}

Following "vector ideas"  \cite{MN_MP_09} we propose as the first step  to write any letter $u_i$ of a plaintext as $n$-tuple ($n$-vector) and after that to apply Procedure \ref{ALG3}. For example it is possible to use a binary representation of characters of the alphabet $A$.

It is possible to divide plain text $u_1, \dots, u_n$ on  parts and to use Procedure \ref{ALG3} to some parts, to a text a part of which has been ciphered by Procedure \ref{ALG3} on a previous ciphering round.

It is possible to change in  Procedure \ref{ALG3}  variables  $x_1, \dots, x_k$ $(1\leq k \leq (n-1))$ by some fixed elements of the set $Q$ and name these elements as leaders. Remark if  $k=n-1$, then we obtain  $n$ chipering images from  any plain-text letter $u$.

In any case the application of only one step Procedure \ref{ALG3} is not very safe since this procedure is not resistant relative to chosen ciphertext attack and chosen plaintext attack.

If in  a system of orthogonal $n$-ary operations there is at least one $n$-ary quasigroup, then we can apply by ciphering of information Algorithm \ref{ALGn_ar} and Procedure \ref{ALG3} together with some non-periodical frequency, i.e., for example, we can apply four times Algorithm \ref{ALGn_ar} and after this we can apply five times  Procedure \ref{ALG3} and so on.

It is possible to  use as a period sequence decimal representation of an irrational or transcendent number.  In this case we can take as a key the sequence of application of Algorithm \ref{ALGn_ar} and Procedure \ref{ALG3}.

 Proposed modifications  make  realization of chosen plaintext attack and chosen ciphertext attack more complicate.

Taking into consideration that in binary case one application of Procedure \ref{ALG3} generates from one plaintext symbol $u$  two cipher symbols, say $v_1, v_2$, we may propose apply Procedure \ref{ALG3} for two plaintext symbols (or to one cipher symbol and one plain symbol, else to two cipher symbols) simultaneously.

The modifications proposed in this subsection need additional researches.

\subsection{Stream-cipher on base of orthogonal system of binary parastrophic quasigroups}

This subsection is more of theoretical than cryptographical character. We propose to use by construction of Algorithm \ref{ALGn_ar} and Procedure \ref{ALG3} orthogonal system of binary parastrophic quasigroups.

We start from the following theorem \cite{MS05}. Here expression
$A\bot {{}^{(23)}A}$ means that  quasigroups $(Q, A)$ and $(Q, {{}^{(23)}A})$ are orthogonal.

\begin{theorem} \label{TH2} For a finite  quasigroup $(Q, A)$ the following equivalences are
fulfilled:

(i)    $ A\bot {{}^{(12)}A}\Longleftrightarrow ((x \backslash  z) \cdot x  = (y \backslash  z)
\cdot y \Longrightarrow x=y)$;

(ii)    $A\bot {{}^{(13)}A} \Longleftrightarrow (zx\cdot x = zy\cdot y \Longrightarrow x=y)
$;

(iii) $ A \bot {{}^{(23)}A} \Longleftrightarrow (x\cdot xz = y\cdot yz \Longrightarrow x=y)
$;

(iv)   $A \bot {{}^{(123)}A} \Longleftrightarrow (x\cdot zx  = y\cdot zy   \Longrightarrow x=y)$;

(v)    $A  \bot {{}^{(132)}A}  \Longleftrightarrow (x z\cdot x
  = yz\cdot y   \Longrightarrow x=y)$

 for all $x, y, z \in Q$.
\end{theorem}

In order to construct quasigroups mentioned in Theorem \ref{TH2}
probably computer search is preferable. It is possible to use   GAP and   Prover \cite{MAC_CUNE_PROV}.

\begin{definition}
A  $T$-quasigroup $(Q,A)$ is a quasigroup of the form $x\cdot y = \varphi x + \psi y + a$, where  $(Q, +)$ is an abelian group, $\varphi, \psi$ are some fixed automorphisms of this group, $a$ is a fixed element of the set $Q$ \cite{pntk, tkpn}.
\end{definition}

In order to construct a quasigroup $(Q, A)$ that is orthogonal with its parastrophe in more theoretical way it is possible to use the following theorem \cite{MS05}.

\begin{theorem} \label{T1}
For a $T$-quasigroup $(Q,A)$ of the form $A(x, y) =  \varphi x + \psi y + a$ over an abelian group $(Q, +)$ the
following equivalences are fulfilled:

(i) $A \bot {}^{12}A \Longleftrightarrow (\varphi - \psi), (\varphi + \psi)$ are permutations of the set $Q$;

(ii) $A\bot {}^{13}A \Longleftrightarrow (\varepsilon + \varphi)$ is a permutation of the set $Q$;

(iii) $A\bot {}^{23}A \Longleftrightarrow (\varepsilon + \psi)$ is a permutation of the set $Q$;

(iv) $A\bot {}^{123}A \Longleftrightarrow (\varphi+ \psi^2)$ is a permutation of the set $Q$;

(v) $A\bot {}^{132}A \Longleftrightarrow (\varphi^2 + \psi)$ is a permutation of the set $Q$.
\end{theorem}

\begin{example}\label{PAR_ORTH_CYCL}
We take the cyclic group $Z_p$, where $p$ is prime. For example, $p=257$ since the number $257$ is prime. Then
$T$-quasigroup $(Q, \circ)$ of the form $x\circ y = k \cdot x + m \cdot y + a$, $k, m, a \in Z_p$, $k, m, k+m, k-m,  k+1, m+1, k^2+m, k+m^2 \neq 0 \pmod{p}$,  where the operation $\cdot$ is multiplication modulo  $p$, is orthogonal to any of its parastrophes.
\end{example}

Any quasigroup $(Q, \circ)$ from Example \ref{PAR_ORTH_CYCL} is orthogonal with  any of its parastrophes. Therefore these quasigroups  are  suitable objects to construct Procedure \ref{ALG3} and, of course, Algorithm \ref{ALGn_ar}.

\noindent \footnotesize
{$^1$Department of Algebra and Number Theory \hfill  $^2$Institute of Mathematics and Computer Science
\\
 E\"otv\"os Lorand University (ELTE) \hfill Academy of Sciences of Moldova  \\
Pazmany Peter setany 1/c, H-1117  Budapest   \hfill  Academiei str. 5, Chi\c{s}in\u{a}u  MD$-$2028 \\
Hungary  \hfill  Moldova  \\
E-mail: \emph{ska@cs.elte.hu \hfill scerb@math.md }}

\end{document}